\documentclass{siamart171218}

\usepackage{lipsum}
\usepackage{amsfonts}
\usepackage{graphicx}
\usepackage{epstopdf}
\usepackage{algorithmic}
\usepackage{epsfig,epsf,xcolor,verbatim}
\usepackage{amsmath,amssymb,float,tikz}

\ifpdf
  \DeclareGraphicsExtensions{.eps,.pdf,.png,.jpg}
\else
  \DeclareGraphicsExtensions{.eps}
\fi

\def\R{{\mathbb R}}
\def\TT{{\mathbb T}}
\def\S{{\mathbb S}}
\def\C{{\mathbb C}}

\def\e{\emptyset}

\def\sd{\S^{d-1}}

\def\ssd{{\sd}\times {\sd}}

\def\dk{\partial K}
\def\dl{\partial L}

\def\tt{\tilde{t}}

\def\rr{{\cal R}}

\def\ot{(\omega,\theta)}

\def\pr{{\rm pr}}

\def\fk{{\cal F}^{(K)}}
\def\fl{{\cal F}^{(L)}}

\def\trapk{\mbox{\rm Trap}(\Omega_K)}
\def\trapl{\mbox{\rm Trap}(\Omega_L)}

\def\sl{{\cal SL}}

\def\endofproof{{\rule{6pt}{6pt}}}

\def\gk{\gamma_K}

\def\la{\left\langle}
\def\ra{\right\rangle}

\def\ep{\epsilon}

\def\ep{\epsilon}

\def\dist{\mbox{\rm dist}}

\def\bs{\bigskip}
\def\ms{\medskip}

\def\bS{{\S}}

\def\trapf{{\footnotesize\rm Trap}}

\def\tF{\widetilde{F}}

\def\tSigma{\widetilde{\Sigma}}

\newsiamremark{remark}{Remark}
\newsiamremark{hypothesis}{Hypothesis}
\crefname{hypothesis}{Hypothesis}{Hypotheses}
\newsiamthm{claim}{Claim}

\headers{Lens Rigidity in Scattering in $\R^2$}{L. Noakes and L. Stoyanov}

\title{Lens Rigidity in Scattering by Unions of Strictly Convex Bodies in $\R^2$\thanks{Submitted to the editors DATE.}}

\author{Lyle Noakes\thanks{Department of Mathematics, University of Western Australia, Crawley WA 6009,
Australia (\email{lyle.noakes@uwa.edu.au}.}
\and Luchezar Stoyanov\thanks{Department of Mathematics, University of Western Australia, Crawley WA 6009,
Australia (\email{luchezar.stoyanov@uwa.edu.au}).}}

\usepackage{amsopn}

\ifpdf
\hypersetup{
  pdftitle={Lens Rigidity in Scattering in $\R^2$},
  pdfauthor={L. Noakes and L. Stoyanov}
}
\fi

\begin{document}

\maketitle

\begin{abstract}
It was proved in \cite{NS1} that obstacles $K$ in $\R^d$  that are finite disjoint 
unions of strictly convex domains with $C^3$ boundaries are uniquely determined by the travelling times 
of billiard trajectories in their exteriors and also by their so called scattering length spectra. However the
case $d = 2$ is not covered in \cite{NS1}. In the present paper we give a separate different proof
of this result in the case $d = 2$.
\end{abstract}



\section{Introduction}

In scattering by an obstacle in $\R^d$ ($d\geq 2$) the obstacle $K$ is a compact subset of ${\R}^d$ 
with a $C^{3}$ boundary $\partial K$ such that  $\Omega_K = \overline{{\R}^d\setminus K}$ is connected. 
A {\it scattering ray} in  $\Omega_K$ is an unbounded in both directions generalized geodesic 
(in the sense of Melrose and  Sj\"ostrand \cite{MS1}, \cite{MS2}). Most of these scattering rays are billiard trajectories 
with finitely many reflection points at $\dk$. In this paper we consider the case when $K$ has the form
\begin{equation} \label{eq:eq1.1}
K = K_1 \cup K_2 \cup \ldots \cup K_{k_0} ,
\end{equation}
where $K_i$ are strictly convex disjoint domains in $\R^d$ with $C^3$ smooth
boundaries $\dk_i$. Then all scattering rays in $\Omega_K$ are billiard trajectories, and the
so called generalized Hamiltonian (or bicharacteristic) flow
$$\fk_t : \S^*(\Omega_K) \longrightarrow \S^*(\Omega_K)$$
coincides with the billiard flow (see \cite{CFS}).

Given an obstacle $K$ in $\R^d$, consider a large ball $M$  containing $K$ in its interior, and let $S_0 = \partial M$
be its boundary sphere.  For any $q\in \dk$ let $\nu_K(q)$ the {\it outward unit normal} to $\dk$.
For $q \in S_0$ we will denote by $\nu(q)$ the {\it inward unit normal} to $S_0$ at $q$. Set
$$\bS^*_+(S_0 ) = \{ x = (q,v) : q \in S_0 \:, \: v \in \sd\: ,\:  \la v, \nu(q) \ra \geq 0\} .$$
Given $x \in  \bS^*_+(S_0)$, define the {\it travelling time} $t_K(x)\geq 0$ as the maximal number (or $\infty$) 
such that $\pr_1(\fk_t(x))$ is in the interior of $\Omega_K \cap M$ for all $0 < t < t_K(x)$, where $\pr_1(p,w) = p$
(see Figure 1).  For $x = (q,v) \in \bS^*_+(S_0)$ with $\la \nu(q),v\ra = 0$ set $t(x) = 0$. The set
$$\{(x;t_K(x)) : x \in \bS^*_+(S_0)\}$$ 
will be called the {\it travelling times spectrum} of $K$.

It is natural to ask what information about the obstacle $K$ can be derived from its travelling times spectrum.
For example: what is the relationship between two obstacles $K$ and $L$ in $\R^d$  if they have (almost) the same
travelling times spectra? We say that $K$ and $L$ {\it have almost the same travelling times} if there exists a subset $R$ of 
full Lebesgue measure in $\bS^*_+(\bS_0)$ such that $t_K(x) =t_L(x)$ for all $x \in R$.

Similar questions can be asked about the so called scattering length spectrum (SLS) of an obstacle. 
Given a scattering ray $\gamma$ in $\Omega_K$, if $\omega\in \sd$ is the incoming direction of $\gamma$ and $\theta\in \sd$ its outgoing direction, 
$\gamma$  will be called an  {\it $\ot$-ray}. For any vector $\xi\in \sd$ denote by $Z_{\xi}$  
{ the hyperplane in ${\R}^d$ orthogonal to $\xi$ and tangent to $S_0$} such that $S_0$ is 
contained in the open half-space $R_{\xi}$ determined by  $Z_{\xi}$ and having $\xi$ as an inner  normal.  
For an $\ot$-ray $\gamma$ in $\Omega$, the {\it sojourn time} $T_{\gamma}$ 
of $\gamma$ is defined by $T_{\gamma} = T'_{\gamma} - 2a$, where $T'_{\gamma}$ is the length of 
the part of $\gamma$ which is contained in $R_{\omega}\cap R_{-\theta}$ and $a$ is the radius of 
$S_0$. It is known that this definition does not depend on the choice of  the sphere $S_0$.  
The {\it scattering length spectrum} of $K$ is defined to be the family  of sets of real numbers 
$$\sl_K =  \{ \sl_K\ot\}_{\ot}$$ 
where $\ot$ runs over $\ssd$ and
$\sl_K\ot$ is the set of sojourn times $T_\gamma$ of all $\ot$-rays $\gamma$ in $\Omega_K$. 
It is known  (cf. \cite{PS}) that for $d \geq 3$, $d$ odd, and $C^\infty$ boundary $\dk$, we have 
$$\sl_K\ot = \mbox{sing supp } s_K(t, \theta, \omega)$$
for almost all $\ot$. Here $s_K$ is the {\it scattering kernel} related to the scattering operator for 
the wave equation in $\R\times \Omega_K$ with Dirichlet boundary condition on 
$\R\times \partial \Omega_K$ (cf. \cite{LP}, \cite{M}, \cite{PS}). Following \cite{St3}, 
we will say that two obstacles  $K$ and $L$ have {\it almost the same SLS} if there exists a subset 
$\rr$ of full Lebesgue measure in $\ssd$  such that  $\sl_K\ot = \sl_L\ot$ for all $\ot\in \rr$.

\bs

\begin{tikzpicture}[xscale=1.44,yscale=0.50] 
  \draw (0,0)  (1,0); 
\draw[blue][thick] (6,-1) arc (0:360:4cm and 7cm); 
\draw[thick] (5,-1) arc (0:360:0.5cm and 0.5cm); 
\draw[thick] (4.1,-2) arc (0:360:0.5cm and 1cm); 
\draw[thick] (5,2) arc (0:360:0.3cm and 0.4cm); 
\draw[thick] (3,2) arc (0:360:0.2cm and 0.1cm); 
\draw[thick] (1,-1) arc (0:360:0.5cm and 0.6cm); 
\draw[thick] (2.1,-2) arc (0:360:0.3cm and 0.6cm); 
\draw[thick] (4,3) arc (0:360:0.2cm and 0.4cm); 
\draw[thick] (2,2) arc (0:360:0.4cm and 0.3cm); 

\draw[thick] (5,-3) arc (0:360:0.4cm and 0.8cm); 
\draw[thick] (4.1,-4) arc (0:360:0.3cm and 0.6cm); 
\draw[thick] (3,-3) arc (0:360:0.3cm and 0.2cm); 
\draw[thick] (3.2,3) arc (0:360:0.3cm and 0.6cm); 
\draw[thick] (1.8,0.4) arc (0:360:0.3cm and 0.8cm); 
\draw[thick] (2.4,0.9) arc (0:360:0.2cm and 0.6cm); 
\draw[thick] (4.2,2) arc (0:360:0.2cm and 0.4cm); 
\draw[thick] (2,-1) arc (0:360:0.1cm and 0.2cm); 

\draw[thick] (-0.5,-1) arc (0:360:0.4cm and 0.5cm); 
\draw[thick] (3.8,1.1) arc (0:360:0.2cm and 0.4cm); 
\draw[thick] (0,-2) arc (0:360:0.3cm and 0.4cm); 
\draw[thick] (-0.2,-3) arc (0:360:0.2cm and 0.4cm); 
\draw[thick] (1,-3.2) arc (0:360:0.5cm and 0.6cm); 
\draw[thick] (2.2,-4) arc (0:360:0.3cm and 0.5cm); 
\draw[thick] (3,-4.5) arc (0:360:0.2cm and 0.4cm); 
\draw[thick] (2,4) arc (0:360:0.4cm and 0.3cm); 

\draw[thick] (5,1) arc (0:360:0.2cm and 0.4cm); 
\draw[thick] (5,0) arc (0:360:0.3cm and 0.2cm); 
\draw[thick] (4.4,0.9) arc (0:360:0.2cm and 0.4cm); 
\draw[thick] (3,1) arc (0:360:0.2cm and 0.2cm); 
\draw[thick] (3.4,0.2) arc (0:360:0.3cm and 0.4cm); 
\draw[thick] (4,-0.2) arc (0:360:0.2cm and 0.6cm); 
\draw[thick] (3,-1) arc (0:360:0.2cm and 0.4cm); 
\draw[thick] (3,-6) arc (0:360:0.2cm and 0.3cm); 

\draw[thick] (2,-6) arc (0:360:0.3cm and 0.5cm); 
\draw[thick] (1.2,-5) arc (0:360:0.3cm and 0.6cm); 
\draw[thick] (1,-6) arc (0:360:0.3cm and 0.2cm); 
\draw[thick] (0.2,-5) arc (0:360:0.3cm and 0.4cm); 
\draw[thick] (-0.6,-4.4) arc (0:360:0.3cm and 0.5cm); 
\draw[thick] (1.2,-1.9) arc (0:360:0.2cm and 0.4cm); 
\draw[thick] (1,2) arc (0:360:0.2cm and 0.4cm); 
\draw[thick] (0.5,3) arc (0:360:0.3cm and 0.2cm); 

\draw[thick] (-0.5,2) arc (0:360:0.4cm and 0.5cm); 
\draw[thick] (-0.2,1.1) arc (0:360:0.2cm and 0.4cm); 
\draw[thick] (0.6,1.2) arc (0:360:0.3cm and 0.7cm); 
\draw[thick] (2.5,3) arc (0:360:0.2cm and 0.4cm); 
\draw[thick] (2,3) arc (0:360:0.2cm and 0.3cm); 
\draw[thick] (1.3,3.2) arc (0:360:0.2cm and 0.3cm); 
\draw[thick] (3.3,4.5) arc (0:360:0.2cm and 0.4cm); 
\draw[thick] (2,2) arc (0:360:0.4cm and 0.3cm);

 \draw (5,4) node[anchor=south west] {$S_0$}; 
\draw (-2.3,-0.4) node[anchor=north west] {$q$};
\draw [->][purple][thick](-1.99,-0.5)--(-1.5,-0.26) ;  \node at (-1.5,0.1) {$v$};
\draw[fill=black] (-1.99,-0.5) ellipse (0.02 and 0.06);
\draw [->][red][thick](-1.99,-0.5)-- (0.3,0.5) ;
\draw[fill=black] (0.3,0.5) ellipse (0.02 and 0.06);
\draw [->][red][thick](0.3,0.5) -- (1.2,0.2) ;
\draw[fill=black] (1.2,0.2) ellipse (0.02 and 0.06);
\draw [->][red][thick](1.2,0.2) -- (0.98,-0.77) ;
\draw[fill=black] (0.98,-0.77) ellipse (0.02 and 0.06);
\draw [->][red][thick](0.98,-0.77) -- (1.67,-1.44) ;
\draw[->][fill=black] (1.67,-1.44)  ellipse (0.02 and 0.06);
\draw [->][red][thick] (1.67,-1.44) --  (2,0.8) ;
\draw[fill=black] (2,0.8) ellipse (0.02 and 0.06);
\draw [->][red][thick] (2,0.8) --  (1.7,1.68) ;
\draw[fill=black] (1.7,1.68) ellipse (0.02 and 0.06);
\draw [->][red][thick] (1.7,1.68) --  (1.45,1.2) ;
\draw[fill=black] (1.45,1.2) ellipse (0.02 and 0.06);
\draw [->][red][thick] (1.45,1.2) --  (1,1.9) ;
\draw[fill=black] (1,1.9) ellipse (0.02 and 0.06);
\draw [->][red][thick] (1,1.9) --  (1.6,3.65) ;
\draw[fill=black] (1.6,3.65) ellipse (0.02 and 0.06);
\draw [->][red][thick] (1.6,3.65) --  (1.8,3.35) ;
\draw[fill=black] (1.8,3.35) ellipse (0.02 and 0.06);
\draw [->][red][thick] (1.8,3.35) --  (2.82,5.83) ;
\draw[fill=black] (2.82,5.83) ellipse (0.02 and 0.06);

\end{tikzpicture}

\bigskip

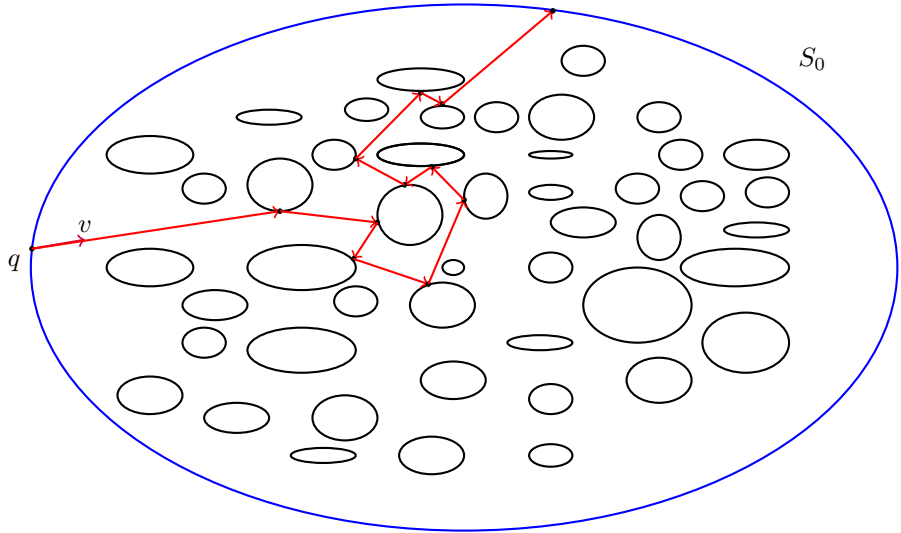
\centerline{Figure 1}

\bigskip

It is a natural and rather important problem in inverse scattering by obstacles to get information 
about the obstacle $K$ from its SLS. It is known that various kinds of information about $K$
can be recovered from its SLS, and for some classes of obstacles $K$ is completely recoverable 
(see \cite{St3}  for more information) -- for example star-shaped obstacles are in this class.

Similar inverse problems concerning metric rigidity have been studied for a very long time in Riemannian geometry 
-- see \cite{SU}, \cite{SUV} and the references there for more information. It appears that some of the methods used in this area, e.g.
those in \cite{Gu}, \cite{DGu}, could be applied to obstacle scattering as well.

More recently various results have been established concerning inverse scattering by obstacles -- see \cite{St2}, \cite{St3},
\cite{NS1} - \cite{NS3}, \cite{St5}. It turns out that some kind of obstacles are uniquely recoverable from their
travelling times spectra and also from their scattering length spectra. For example, it was shown in \cite{NS1} that if $K$ and 
$L$ are finite disjoint unions of strictly convex bodies in $\R^d$ with $C^3$ boundaries and $K$ and $L$ have almost the same 
travelling times spectra (or almost the same SLS), then $K = L$. However the argument in \cite{NS1} does not work in the case $d = 2$. 
We are grateful to Antoine Gansemer who pointed this to us. As he showed in \cite{Gan}, when $d = 2$ and $k_0 > 1$ the set 
$\bS^*_+(S_0)\setminus \trapk$ is disconnected, and then the argument in \cite{NS1} does not work. Here  $\trapk$ is the 
{\it set of all trapped points} in $\S^*(\Omega_K)$, i.e. points $x = (q,v)\in  \bS^*(\Omega_K)$ such that either the forward billiard trajectory 
$$\gk^+(x) = \{\pr_1(\fk_t(x)) : t \geq 0\} $$ 
or the backward trajectory $\gk^- (q,v) = \gk^+(q, -v)$ is infinitely long.
That is, either the billiard 
trajectory in the exterior of $K$ issued from $q$ in the direction of $v$ is bounded (contained entirely in $M$) or the 
one issued from $q$ in the direction of $-v$ is bounded. The obstacle $K$ is called {\it non-trapping} if $\trapk = \e$.

Here we prove the following.

\bs

\begin{theorem}\label{thm:Thm 1.1}
Let $K$ and $L$ be obstacles  in $\R^2$ such that each of them is a  finite disjoint union of strictly convex compact domains  
with $C^3$ boundaries. If $K$ and $L$ have  almost the same travelling times or almost the same scattering length spectra, 
then $K = L$.
\end{theorem}

The argument we use is completely different from that in \cite{NS1}.
A result similar to that in \cite{NS1} concerning non-trapping
obstacles satisfying certain non-degeneracy conditions  was proved recently in \cite{St5}.

The set of trapped points plays a rather important role in various inverse problems in scattering by obstacles, and also
in problems on metric rigidity in Riemannian geometry. It is known that $\trapk \cap \bS^*_+(S_0)$ has Lebesgue
measure zero in $\bS^*_+(S_0)$ (see Sect. 4 for more information about this). However, as an example of M. Livshits shows  
(see  Ch. 5 in \cite{M} or Figure 1  in \cite{NS1}),  in general the set of points $x\in \bS^*(\Omega_K)$ for which 
$$\gk(x) = \gk^+(x) \cup \gk^-(x) $$
is  trapped in both directions   may contain a non-trivial open set. In the latter case the  obstacle cannot be recovered from travelling 
times (and also from the SLS). Similar examples in higher dimensions are given  in \cite{NS3}.

\begin{definition}
Let $K, L$ be two obstacles  in $\R^d$.
We will say that $\Omega_K$ and $\Omega_L$  {\it have conjugate flows} if  there exists a homeomorphism 
$$\Phi : \S^*(\Omega_{K})\setminus \trapk  \longrightarrow  \S^*(\Omega_{L})\setminus\trapl$$
which is $C^1$ on an open dense subset of  $\S^*(\Omega_{K})\setminus \trapk$ and  satisfies 
$$\fl_t\circ \Phi = \Phi\circ \fk_t \quad, \quad t\in \R ,$$
and $\Phi = \mbox{id}$ on  $\S^*(\R^d\setminus M)\setminus \trapk = \S^*(\R^d\setminus M)\setminus\trapl$.
\end{definition}

For $K, L$ in a generic class of obstacles in $\R^d$ ($d\geq 2$), which includes the type of obstacles considered here,
it is known  that if $K$ and $L$ have 
almost the same SLS or almost the same travelling times, then $\Omega_K$ and $\Omega_L$ have conjugate 
flows (\cite{St3} and \cite{NS2}). Thus, \cref{thm:Thm 1.1} is an immediate consequence of the following.

\ms

\begin{theorem}\label{thm:Thm 1.3}
Let  each of the obstacles $K$ and $L$ be a 
finite disjoint union of strictly convex compact domains in $\R^2$ with $C^3$ boundaries. 
If $\Omega_K$ and  $\Omega_L$  have conjugate flows, then $K = L$.
\end{theorem}

\ms

We prove \cref{thm:Thm 1.3} in Sect. 3 below. In Sect. 2 we state some useful results from \cite{St2} and \cite{St3}.
It turns out that billiard trajectories with tangent points to the boundary of the obstacle play an important role
in the two-dimensional case considered here. 
In Sect. 4 we prove that the set of trapped points $\trapk$ has Lebesgue measure zero in $\S^*(\Omega_K)$.

\section{Preliminaries}

Next, we describe some propositions from \cite{St2} and \cite{St3} that are needed in the proof of \cref{thm:Thm 1.3}.
We state them in the general case $d \geq 2$, although later on we will use them in the special case  $d = 2$.

\ms

\noindent 
{\bf Standing Assumption.} $K$ and $L$ are finite disjoint unions of strictly convex domains
in $\R^d$ ($d\geq 2$)  with $C^3$ boundaries and with conjugate flows $\fk_t$ and $\fl_t$.
 
\ms 

\begin{proposition} \label{prop:Prop 2.1}
(\cite{St2})
(a) {\it There exists a countable family $\{ M_i\} = \{ M_i^{(K)}\}$
of codimension $1$ submanifolds of $\S^*_+(S_0)\setminus \trapk$ such that every 
$$\sigma \in \S^*_+(S_0)\setminus (\trapk \cup_i M_i)$$ 
generates a simply reflecting ray in $\Omega_K$. Moreover the family  $\{ M_i\}$
is locally finite, that is any compact subset of  $\S^*_+(S_0)\setminus \trapk$ has common points with
only finitely many of the submanifolds $M_i$.} 

\ms

(b) {\it There exists a countable family $\{ R_i\}$
of codimension $2$ smooth submanifolds of $\S^*_+(S_0)$ such that for any
$\sigma\in \S^*_+(S_0)\setminus (\cup_i R_i)$ the trajectory $\gk(\sigma)$ has at most one
tangency to $\dk$.}

\ms

(c)  {\it There exists a countable family $\{ Q_i\}$ of codimension $2$ smooth submanifolds of 
$\S^*_{\dk} (\Omega_K)$ such that for any $\sigma\in \S^*_+(\dk)\setminus (\cup_i Q_i)$ the trajectory 
$\gk(\sigma)$ has at most one tangency to $\dk$.}
\end{proposition}

\ms

It follows from the conjugacy of flows and Proposition 4.3 in \cite{St3} that
the submanifolds $M_i$ are the same for $K$ and $L$, i.e. $M_i^{(K)} = M_i^{(L)}$ for all $i$.

The following is Lemma 5.2 in \cite{St2}. In fact the lemma in  \cite{St2} 
assumes $C^\infty$ smoothness for the submanifold $X$, however its proof only
requires $C^3$ smoothness. 

\bigskip

\begin{proposition}  \label{prop:Prop 2.2}
 Let $X$ be a $C^3$ smooth submanifold of codimension $1$ in $\R^d$ , and let $x_0\in X$ and 
$\xi_0\in T_xX$, $\| \xi_0\| = 1$, be such that the normal curvature of $X$ at  $x_0$ in the direction $\xi_0$ is non-zero.
Then for every $\epsilon > 0$ there exist an open neighbourhood $V$ of  $x_0$ in $X$, a smooth map 
$$V \ni x \mapsto \xi(x) \in T_x X$$ 
and a smooth positive  function $t(x) \in [\delta,\epsilon]$ on $V$ for some $\delta > 0$ such that
$$Y = \{y(x) = x + t(x)\xi(x) : x\in V\}$$ 
is a smooth strictly convex surface with an unit normal field $\nu_Y(y(x)) = \xi(x)$, $x\in V$.
That is, the normal field of $Y$ consists of vectors tangent to $X$  at the corresponding points of $V$. (See Figure 2.)
\end{proposition}

\bs

As one would expect, the case $d = 2$ of the above proposition is rather easy to prove.

\bs

\begin{tikzpicture}[xscale=1.44,yscale=0.50] 
  \draw (0,0)  (1,0); 
  \draw[blue][thick] (6,-1) arc (0:180:2cm and 2cm); 
  \draw[purple][thick] (6,0) arc (-30:30:0.5cm and 2cm); 
  \draw (2,0) node[anchor=south west] {$X$}; 
  \draw (3.85,0.75) node[anchor=north west] {$x_0$};
  \draw (6,0.3) node[anchor=north west] {$Y$};
\draw [->](3,1)--(7,1) ;  \node at (7.2,1) {$\xi_0$};
\draw [->](3,0.9)--(7,1.5) ;  
\draw [->](3,0.8)--(7,1.9) ;  
\draw [->](3,1.2)--(7,0.5) ;  
\draw[fill=black] (4,1) ellipse (0.02 and 0.06);

\end{tikzpicture}

\bigskip

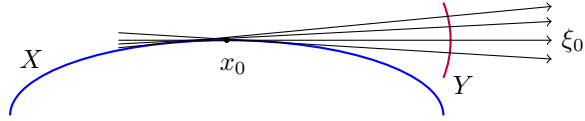
\centerline{Figure 2}

\bs

An important consequence of the above is the following proposition which can be proved using part of the argument in the proof of 
Proposition 5.5 in \cite{St2}. For completeness we sketch the proof in the Appendix.

\bs

\begin{proposition}  \label{prop:Prop 2.3}
Let $K$ be an obstacle in $\R^2$ which is a  finite disjoint union of strictly convex compact domains  
with $C^3$ boundaries. Then 
$$\dim(\S^*(\dk) \cap \trapk) = 0 .$$
In particular, $\S^*(\dk) \cap \trapk$ does not contain non-trivial open subsets of $\S^*(\dk)$.
\end{proposition}

\ms

Here we denote by $\dim(X)$ the {\it topological dimension} of a subset $X$ of $\R^2$.

It turns out that for the type of obstacles considered in this paper  the set $\trapk$ of trapped points has 
Lebesgue measure zero in $\S^*(\Omega_K)$. While formally this fact is not necessary for the proof of  \cref{thm:Thm 1.3}, 
we mention it here since it is a rather important feature of the billiard flow in the case considered in this paper (and also
in \cite{NS1}, \cite{St3}, etc.).
This  appears to be accepted as a `known fact' although we could not find a formal proof anywhere in the literature. However a simple 
proof follows from known facts, e.g. using the ergodicity of the
so called dispersive (Sinai) billiards (see \cite{Si1}, \cite{Si1}).  

\ms

\begin{proposition}\label{prop:Prop 2.4}
Let $K$ be an obstacle in $\R^d$ of the form \cref{eq:eq1.1}. Then  the set $\trapk$ of all of trapped points of $\S^*(\Omega_K)$ has 
Lebesgue measure zero in $\S^*(\Omega_K)$.
\end{proposition}

\ms

We  provide a proof of this proposition in Sect. 4 below.






\section{Proof of  \cref{thm:Thm 1.3}}
\setcounter{equation}{0}

Assume that the obstacles $K$ and $L$ in $\R^2$ satisfy the assumptions of  \cref{thm:Thm 1.3}. 

We claim that $K \subset L$. Assume this is not true and fix an arbitrary $x_0 \in \dk$ such that 
$x_0 \notin L$. Let $\xi_0 \in \S^1$ be one of the unit vectors tangent to $\dk$ at $x_0$. 

 It follows from \cref{prop:Prop 2.2} that there exists a small $\ep_0 > 0$,
an open neighbourhood $V_0$ of  $x_0$ in $\dk$, a $C^2$ map  $V_0 \ni x \mapsto \xi(x) \in \S^*_x (\dk)$ 
and a $C^2$  positive function $t(x) \in [\delta,\ep_0]$ on $V_0$ for some $\delta \in (0,\ep_0)$ such that
$$\Sigma = \{y(x) = x + t(x)\xi(x) : x\in V_0\}$$ 
is a $C^2$ strictly convex curve with unit normal field $\nu_\Sigma (y(x)) = \xi(x)$, $x\in V_0$. So, for 
any $x\in V_0$ the straight line through $y(x)$ with direction $\xi(x)$ is tangent to $\dk$ at $x$.
Set $y_0 = x_0 + \ep_0 \xi_0 \in \Sigma$. 

It follows from \cref{prop:Prop 2.3} that for the subset
$$\Sigma' = \{ y\in \Sigma : (y, \nu_\Sigma(y)) \notin \trapk\}$$
we have $\dim(\Sigma \setminus \Sigma') = 0$. Thus, $\dim(\Sigma') = 1$.

Next, \cref{prop:Prop 2.1}  implies that for all but countably many 
$y \in \Sigma'$ the trajectories $\gamma_K(y, \nu_\Sigma(y))$ and $\gamma_L(y, \nu_\Sigma(y))$
have at most one tangency to $\dk$ and $\dl$, respectively. For such $y$, since 
$\gamma_K(y, \nu_\Sigma(y))$ has a tangent point to $\dk$, it must have exactly one tangent point to 
$\dk$. Since the flows $\fk_t$ and $\fl_t$ are conjugate by assumption,  $\gamma_L(y, \nu_\Sigma(y))$ also 
must have exactly one tangent point  $z(y)$ to $\dl$. More precisely, if $(y, \nu_\Sigma(y)) = \fk_t(\sigma)$
for some $\sigma \in \S^*_+(S_0)$ and some $t > 0$, then the travelling time function $t_K$ has 
a singularity at $\sigma$. Since $t_K = t_L$ on $\S^*_+(S_0)$ near $\sigma$, the function $t_L$ also has
a singularity at $\sigma$, so $\gamma_L(y, \nu_\Sigma(y)) = \gamma_L(\sigma)$ must have a tangent point to $\dl$.


Assume for a moment that for every $z \in \dl$ there exists an open neighbourhood $W_z$  of $z$ in $\dl$ such that 
$$W_z \cap \{ z(y) : y \in \Sigma'\}$$ 
has  topological dimension zero.
Covering $\dl$ with a finite number of neighbourhoods $W_z$, it follows that $\Sigma'$ has topological dimension zero -- a contradiction. 
Thus, there exists $z_0 \in \dl$ such that for every open neighbourhood $W_0$  of $z_0$ in $\dl$ the set
$$W_0 \cap \{ z(y) : y \in \Sigma\}$$ 
has  topological dimension one.  Replacing $y_0$ (and therefore $x_0$ as well) by an appropriate
nearby point on $\Sigma'$, we may assume that $z_0 = \pr_1(\fl_{t_0}(y_0, \nu_\Sigma(y_0)))$ for some  $t_0\in \R$, $t_0\neq 0$.

We will assume that $t_0 > 0$; otherwise we just have to replace $\xi_0$ by $-\xi_0$ and the curve
$\Sigma$ by $\{x - t(x)\xi(x) : x\in V_0\}$. Let  
$$\fl_{t_0}(y_0, \nu_\Sigma(y_0)) = (z_0, -\zeta_0) .$$

\bs

\begin{tikzpicture}[xscale=1.50,yscale=0.50] 
  \draw (0,0)  (1,0); 
  \draw[blue][thick] (2,-1) arc (0:180:1cm and 1cm); 
  \draw[purple][thick] (2.2,-1) arc (-30:30:0.5cm and 2cm); 
  \draw (0.2,-1.2) node[anchor=south west] {$\dk$}; 
  \draw (0.8,-0.2) node[anchor=north west] {$x_0$};
  \draw (2.2,-0.8) node[anchor=north west] {$\Sigma$};
\draw [->](1,0)--(2.6,0) ;  \node at (2.7,-0.25) {$\xi_0$};
\draw[fill=black] (1,0) ellipse (0.02 and 0.06);

\draw[fill=black] (7,0) ellipse (0.02 and 0.06);
\draw[blue][thick] (8,-1) arc (0:180:1cm and 1cm); 
\draw[purple][thick] (5.7,1) arc (150:210:0.5cm and 2cm); 
  \draw (7.4,-1.2) node[anchor=south west] {$\dl$}; 
  \draw (6.8,-0.2) node[anchor=north west] {$z_0$};

 \draw (5.4,-0.9) node[anchor=north west] {$X$};
 \draw (5.7,-0.45) node[anchor=north west] {$p$};

\draw [->](7,0)--(5.2,0) ;  \node at (5.1,-0.25) {$\zeta_0$};
\draw[-] (2,0.5)--(2.8,0.8) ; \draw[fill=black] (2.8,0.8) ellipse (0.02 and 0.06);
\draw[fill=black] (2.25,0.58) ellipse (0.02 and 0.06);
\draw[-] (2.8,0.8)--(3.8,-0.6) ; \draw[fill=black] (3.8,-0.6) ellipse (0.02 and 0.06);
\draw[-] (3.8,-0.6)--(4.2,1.6) ; \draw[fill=black] (4.2,1.6) ellipse (0.02 and 0.06);
\draw[-] (4.2,1.6)--(4.8,-0.9) ; \draw[fill=black] (4.8,-0.9) ellipse (0.02 and 0.06);
\draw[-] (4.8,-0.9)--(5.9,-0.5) ; \draw[fill=black] (5.65,-0.58) ellipse (0.02 and 0.06);

\end{tikzpicture}

\bigskip

\bs

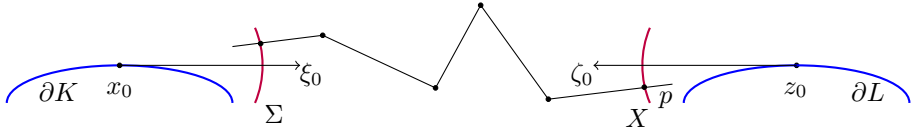
\centerline{Figure 3}

\bigskip

Using again \cref{prop:Prop 2.2}, assuming $\ep_0 > 0$ is sufficiently small and shrinking  the open 
neighbourhood $W_0$ of  $z_0$ in $\dl$ if necessary, there exist a $C^2$ map  
$$W_0 \ni z \mapsto \zeta(z) \in \S^*_z (\dl)$$
and a $C^2$  positive function $s(z) \in [\delta,\ep_0]$ on $W_0$ for some $\delta \in (0,\ep_0)$ such that
$\zeta(z_0) = \zeta_0$ and 
$$X = \{p(z) = z + s(z)\, \zeta(z) : z\in W_0\}$$ 
is a $C^2$ strictly convex curve with unit normal field $\nu_X (p(z)) = \zeta(z)$, $z\in W_0$.   So, for 
any $z\in W_0$ the straight line through $p(z)$ with direction $\zeta(z)$ is tangent to $\dl$ at $z$
(see Figure 3).
Set $p_0 = z_0 + \ep_0 \zeta_0 \in X$.

We will now use a basic property of dispersive (Sinai) billiard flows concerning propagation of convex fronts.
Let 
$$x_1 = \pr_1(\fl_{t_1}(x_0,\xi_0)), \ldots, x_k = \pr_1(\fl_{t_k}(x_0,\xi_0))$$ 
be the common points of $\gamma_L(x_0,\xi_0)$ with $\dl$ (if any) with $0 < t_1 < \ldots < t_k < t_0$. 
Fix an arbitrary $T \in (t_k,t_0)$ close to $t_k$. It then follows from a well-known result of Sinai (\cite{Si1}; see also
\cite{Si2}) that there exists an open neighbourhood $\Sigma_0$ of $y_0$ in $\Sigma$ such that
$$Y = \{ \pr_1(\fl_{T}(y,\nu_\Sigma(y))) : y \in \Sigma_0\}$$
is a strictly convex curve in $\R^2$ with a unit normal field 
$$\nu_Y( y,\nu_\Sigma(y))  = \pr_2(\fl_{T}(y,\nu_\Sigma(y))) .$$

\bs

\begin{tikzpicture}[xscale=1.30,yscale=0.40] 
  \draw (0,0)  (1,0); 
  \draw[blue][thick] (2,-1) arc (0:180:1cm and 1cm); 
  \draw[purple][thick] (2.05,-1) arc (-30:30:0.5cm and 2cm); 
  \draw (0.2,-1.4) node[anchor=south west] {$\dk$}; 
  \draw (0.8,-0.2) node[anchor=north west] {$x_0$};
  \draw (1.85,2.2) node[anchor=north west] {$\Sigma_0$};
\draw [->](1,0)--(1.7,0) ;  \node at (1.6,0.8) {$\xi_0$};
\node at (2.3,-0.6) {$y_0$};
\draw[fill=black] (1,0) ellipse (0.02 and 0.06);
\draw[fill=black] (2.12,0) ellipse (0.02 and 0.06);

\draw[fill=black] (8,0) ellipse (0.02 and 0.06);
\draw[blue][thick] (9,-1) arc (0:180:1cm and 1cm); 
\draw[purple][thick] (6.95,1) arc (150:210:0.5cm and 2cm); 
  \draw (8.5,-1.8) node[anchor=south west] {$\dl$}; 
  \draw (7.8,-0.2) node[anchor=north west] {$z_0$};

 \draw (6.8,2.2) node[anchor=north west] {$X$};
\node at (6.75,-0.6) {$p_0$};
\draw[fill=black] (6.88,0) ellipse (0.02 and 0.06);

\draw [->](8,0)--(7.3,0) ;  \node at (7.5,0.8) {$\zeta_0$};
\draw[-] (1.7,0)--(2.8,0) ; 
\draw[fill=black] (2.8,0) ellipse (0.02 and 0.06);
\node at (2.8,0.8) {$x_1$};
\draw[-] (2.8,0)--(3.8,-1.3) ; 
\draw[fill=black] (3.8,-1.3) ellipse (0.02 and 0.06);
\node at (4,-1.6) {$x_2$};
\draw[-] (3.8,-1.3)--(4.5,2) ; 
\draw[fill=black] (4.5,2) ellipse (0.02 and 0.06);
\draw[-] (4.5,2)--(5.2,0) ; 
\draw[fill=black] (5.2,0) ellipse (0.02 and 0.06);
\draw[-] (5.2,0)--(7.3,0) ; 
\draw[violet][thick] (5.6,-1.7) arc (-60:60:0.3cm and 2cm); 
\draw[fill=black] (5.73,0) ellipse (0.02 and 0.06);
 \draw (5.5,2.85) node[anchor=north west] {$Y$};
 \node at (5.55,-0.6) {$q_0$};
 \draw [->](5.73,0)--(6.43,0) ;  \node at (6.1,0.8) {$\eta_0$};
 \node at (5.2,-0.6) {$x_k$};

\end{tikzpicture}

\bigskip

\bs

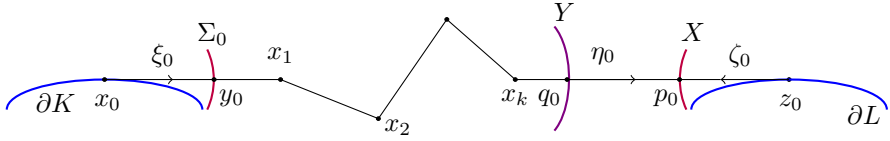
\centerline{Figure 4}

\bigskip

Set 
$$q_0 = \pr_1(\fl_{T}(y_0,\nu_\Sigma(y_0))) \in Y \quad, \quad \eta_0 = \nu_Y( q_0)$$ 
(see Figure 4).
It follows from the constructions 
of $\Sigma$, the point $z_0 \in \dl$, the neighbourhood $W_0$ and the convex fronts $X$ and $Y$ that for 
$y \in \Sigma_0 \cap \Sigma'$  the point 
$$q = \pr_1(\fl_{T}(y,\nu_\Sigma(y))) \in Y$$ 
is such that the straightline ray issued from $q$ in direction $\nu_Y( q)$ hits $X$ perpendicularly. However, due to
the strict convexity of $X$ and $Y$, this is only possible when $y = y_0$; a contradiction.

This proves that we must have $K \subset L$.

Using a similar argument we derive that $L \subset K$, as well. Therefore $K = L$.
\endofproof

\bs

\section{On the set of trapped points}
\setcounter{equation}{0}

Here we prove \cref{prop:Prop 2.4}.

Assume again that $K$ is an obstacle in $\R^d$ ($d\geq 2$) of the form \cref{eq:eq1.1} where $K_i$ are strictly convex disjoint domains in 
$\R^d$ with $C^3$ smooth boundaries $\dk_i$.  Let $\lambda$ be the Lebesgue measure on $\S^*(\R^d)$.  
Let $S_0$ be a large sphere in $\R^d$ as in Sect. 1, and let  $\mu$ be the 
{\it Liouville measure } on $\S^*_+(S_0)$ defined by 
$$d\mu = d\rho(q) d\omega_q |\la \nu(q), v\ra| ,$$
where $\rho$ 
is the measure on $S_0$ determined by the Riemannian metric on $S_0$ and $\omega_q$ is the Lebesgue measure on the
$(d-2)$-dimensional sphere $\S_q(S_0)$ (see e.g. Sect. 6.1 in \cite{CFS}).

We will need the following generalisation of Santalo's formula proved in \cite{St4}. In fact, the latter deals with general billiard flows on Riemannian 
manifolds (under some natural assumptions),  however here we will restrict ourselves to the case considered in Sect. 1.

\ms

\begin{theorem}\label{thm:Thm 4.1}
(\cite{St4})
Let $K$ be as above. Then for every $\lambda$-measurable function 
$$f : \S^*(\Omega_K) \setminus \trapk \longrightarrow \C$$
such that $|f|$ is integrable we have
\begin{eqnarray*}
\displaystyle
     \int_{\S^*(\Omega_K)\setminus \trapf (\Omega_K)} f(x) \;  d\lambda(x) 
     =  \int_{\S^*_+(S_0) \setminus \trapf(\Omega_K)} \left( \int_0^{t_K(x)} f(\fk_t(x)) \;dt \right)\, d\mu(x) .
\end{eqnarray*}
\end{theorem}

\ms

As we mentioned earlier $\trapk \cap \S^*_+(S_0)$ has Lebesgue measure zero in $\S^*_+(S_0)$ (see Theorem 1.6.2 in \cite{LP}; see also
Proposition 2.3 in \cite{St2} for a more rigorous proof).
Using this and the above theorem with $f =1$ gives the following.

\ms

\begin{corollary}\label{cor:Cor 4.2}
(\cite{St4})
Let $K$ be as above. Then 
\begin{eqnarray*}
\displaystyle
     \lambda (\S^*(\Omega_K)\setminus \trapk)
     =  \int_{\S^*_+(S_0) \setminus \trapf(\Omega_K)} t_K(x)\, d\mu(x) .
\end{eqnarray*}
That is, 
\begin{eqnarray*}
\displaystyle
      \lambda (\trapk)  =  \lambda (\S^*(\Omega_K)) -  \int_{\S^*_+(S_0)} t_K(x)\, d\mu(x) .
\end{eqnarray*}

\end{corollary}

\ms

\noindent
{\it Proof of } \cref{prop:Prop 2.4}.
We can regard $K$ as a subset of a domain $Q$ in $\R^d$ with a piecewise smooth boundary which is strictly convex inwards\footnote{Or as a 
domain on the flat $d$-dimensional torus $\TT^d$. Both embeddings will produce the required result.} (see Figure 5). Consider the billiard flow $\phi_t$
on $\S^*(Q)$. It is well-known (see \cite{CFS}) that $\phi_t$ preserves the Lebesgue measure $\lambda$ (restricted to $\S^*(Q)$). 
Moreover $\phi_t$ is ergodic with respect to $\lambda$ (\cite{Si1}, \cite{Si2}).

\bs

\begin{tikzpicture}[xscale=1.44,yscale=0.50] 
  \draw (0,0)  (1,0); 
\draw[blue][thick] (1,-3.7) arc (-60:60:1cm and 5cm); 
\draw[blue][thick] (7.3,6.05) arc (104:250:1cm and 5cm); 
\draw[blue][thick] (1,5) arc (240:320:5cm and 5cm); 
\draw[blue][thick] (7.22,-3.58) arc (60:122.5:6cm and 8cm);

\draw[thick] (6.2,2) arc (0:360:0.5cm and 0.5cm); 
\draw[thick] (5.5,-1) arc (0:360:0.5cm and 1cm); 
\draw[thick] (5,2) arc (0:360:0.3cm and 0.4cm); 

\draw[thick] (3.8,1.1) arc (0:360:0.2cm and 0.4cm); 
\draw[thick] (4,-1) arc (0:360:0.3cm and 0.4cm); 
\draw[thick] (2.8,1) arc (0:360:0.3cm and 0.5cm); 
\draw[thick] (3,3) arc (0:360:0.2cm and 0.4cm); 

 \draw (7,4) node[anchor=south west] {$\partial Q$}; 
  \draw (4,0) node[anchor=south west] {$K$}; 

\end{tikzpicture}

\bigskip

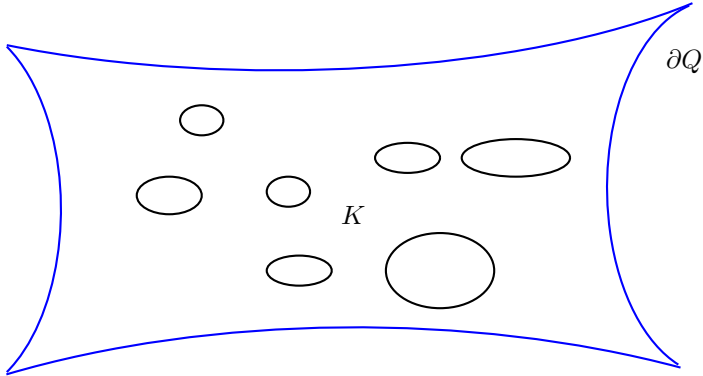
\centerline{Figure 5}

\bigskip

Let $T$ be the set of points $x \in \S^*(\Omega_K)$ such that the billiard trajectory $\gamma_K(x)$ is trapped in both directions.
Then \cref{cor:Cor 4.2}  and the fact mentioned above that $\trapk \cap \S^*_+(S_0)$ has Lebesgue measure zero in $\S^*_+(S_0)$  imply that
$\trapk\setminus T$ has Lebesgue measure zero in $\S^*(\Omega_K)$. So, it is enough to prove that $\lambda(T) = 0$.

The billiard flow $\phi_t$ coincides with the flow $\fk_t$ on the set $T$, and $T$ is an invariant set
with respect to $\fk_t$, and so with respect to $\phi_t$.  Clearly $T$ is a proper subset of $S^*(Q)$ and $S^*(Q) \setminus T$ 
has positive measure. Now the ergodicity of $\phi_t$ implies that $\lambda(T) = 0$.
\endofproof

\bs

\section{Appendix}
\setcounter{equation}{0}

Here we prove \cref{prop:Prop 2.3} using part of the argument in the proof of Proposition 5.5 in \cite{St2}.

It is enough to prove that every $x_0 \in \dk$ has an open neighbourhood $V_0$ in $\dk$ such that
$\dim(\S^*(V_0)\cap \trapk) = 0$.

Let $x_0\in \dk$. As in the proof of \cref{thm:Thm 1.3}, it follows from \cref{prop:Prop 2.2} that there exists a small $\ep_0 > 0$,
an open neighbourhood $V_0$ of  $x_0$ in $\dk$, a $C^2$ map  $V_0 \ni x \mapsto \xi(x) \in \S^*_x (\dk)$ 
and a $C^2$  positive function $t(x) \in [\delta,\ep_0]$ on $V_0$ for some $\delta \in (0,\ep_0)$ such that
$$\Sigma = \{y(x) = x + t(x)\xi(x) : x\in V_0\}$$ 
is a $C^2$ strictly convex curve with unit normal field $\nu_\Sigma (y(x)) = \xi(x)$, $x\in V_0$. Set
$$\tSigma = \{ (y, \nu_\Sigma(y)) : y  \in \Sigma \} .$$

It follows from \cref{prop:Prop 2.1}(c) that there exists a countable subset
$X' = \{Q_i\}$ of  $\S^*(\dk)$ such that for any $\sigma\in S^*(\dk) \setminus X'$, the trajectory $\gk(\sigma)$
has at most one tangency to $\dk$, and therefore it has exactly one tangency to $\dk$. 

Let $X_0$ the set of those $\sigma \in \tSigma\cap \trapk$ such that
the trajectory $\gk^+(\sigma)$ has no tangencies to $\dk$.
Set $F = \{ 1, 2, \ldots, k_0\}$, and consider
$$\tF = \prod_{r=1}^\infty F$$ 
with the product topology. It is well known that $\dim (\tF) = 0$ and therefore every subspace
of $\tF$ has topological dimension zero (cf. e.g. \cite{HW} or \cite{E}).
Consider the map $f : X_0 \longrightarrow \tF$, defined by
$$f(\sigma) = (i_1, i_2, \ldots, i_n, \ldots) ,$$
where the $n$th reflection point
of $\gk^+(\sigma)$ belongs to $\dk_{i_n}$ for all $n = 1,2, \ldots$.
Clearly, the map $f$ is continuous and it follows from \cite{St1} that $f$ is injective, so it
defines a homeomorphism $f : X_0 \longrightarrow f(X_0)$.
Thus, $X_0$ is homeomorphic to a subspace of $\tF$ and therefore $\dim(X_0) = 0$.

Now the Sum  Theorem for $\dim$ (cf. \cite{HW} or \cite{E}) shows that $\dim( X' \cup X_0) = 0 .$ Since 
$S^*(V_0) \cap \trapk$ is naturally homeomorphic to $ X' \cup X_0$,  it follows that 
$$\dim(S^*(V_0) \cap \trapk) = 0 .$$
This proves the proposition. 
\endofproof

\bs


\newpage
\bs

\bibliographystyle{siamplain}

\end{document}